\documentclass[11pt,leqno]{article}
\usepackage{amssymb}
\usepackage{amsmath}
\usepackage{latexsym}
\usepackage{amsfonts}

\begin{document}
%%%%%%%%%%%%%%%%%%%
\newcommand{\rs}{{\Phi}}   		%rs.tex
\newcommand{\ov}{{\cal O}(V)} 		%frg.tex
\renewcommand{\ss}{{\Delta}}		%pss.tex
\newcommand{\ps}{{\Pi}}			%pss.tex
\newcommand{\ess}{{\delta}}		%conj.tex
\newcommand{\len}{{\ell}}		%length.tex
\renewcommand{\mod}[1]{|#1|}		%length.tex
\newcommand{\fg}{{\frak g}}		%intro.tex
\newcommand{\frg}{{\frak g}}		%intro.tex
%%%%%%%%%%%%%%%%%%%%%%%%%%%%%%%%%%%%%%%%%%%%%%%%%%%%%%
%%%%%%%%%%%%%%%%%%%%%%%%%%%%%%%%%%%%%%%%%%%%%%%%%%%%%%
\newcommand{\st}{{\,|\,}}		%general
\renewcommand{\frak}{\mathfrak}		%general
\renewcommand{\Bbb}{\mathbb}		%general
%%%%%%%%%%%%%%%%%%%%%%%%%%%%%%%%%%%%%%%%%%%%%%%%%%%%%%

\newcommand{\state}{{\phi}}		%state.tex
\newcommand{\states}{{\cal{S}}}		%state.tex
\newcommand{\hilb}{{\cal{H}}}		%state.tex
\newcommand{\sleq}{{\leq_1}}		%state.tex
\newcommand{\sgeq}{{\geq_1}}		%state.tex

\renewcommand{\lneq}{{<}}		%global.tex

\newcommand{\shape}{{\lambda}}		%monom.tex
\newcommand{\monom}{M}			%monom.tex
\newcommand{\ul}{\underline}		%theorem.tex
\newcommand{\bc}{{\Bbb{C}}}		%theorem.tex
\newcommand{\hit}{{\cal S}}		%prop.tex
\newcommand{\tnot}{{{\ul{t}}^{0}}}	%prop.tex
\newcommand{\ult}{{\ul{t}}}		%prop.tex
\newcommand{\ultp}{{\ul{t}'}}		%prop.tex

\newcommand{\ds}{\displaystyle}		%prop.tex
\newcommand{\be}{\begin{enumerate}}     %prop.tex
\newcommand{\ee}{\end{enumerate}}	%prop.tex
\newcommand{\weight}{{\mbox{\rm wt}}}	%comb.tex
\newcommand{\uqminus}{{U^-_q(\frak g)}} %intro.tex
\newcommand{\simga}{\sigma}
\newcommand{\lamdba}{\lambda}
\newcommand{\SMT}{ standard monomial theoretic }
\newcommand{\tensor}{\otimes}
\newcommand{\linfty}{{\cal L}(\infty)}	%intro.tex
\newcommand{\ratf}{\Bbb{Q}(q)}		%prop.tex
\newcommand{\vz}{V_{\Bbb{Z}}}		%prop.tex
\newcommand{\Ell}{{\cal L}}		%prop.tex
\newcommand{\latz}{\Ell_{\Bbb{Z}}}	%prop.tex
\newcommand{\vzi}{V_\Bbb{Z}(i)}		%prop.tex
\newcommand{\uzminus}{U_{\Bbb{Z}}^-}	%prop.tex
\newcommand{\surjection}{\twoheadrightarrow}%theorem.tex
\newcommand{\matp}{{\underline{\underline{P}}}}  %theorem.tex
\newcommand{\matb}{{\underline{\underline{B}}}}  %theorem.tex
\newcommand{\matc}{{\underline{\underline{C}}}} %theorem.tex
\newcommand{\matd}{{\underline{\underline{D}}}}%theorem.tex
\newcommand{\maxi}{{\frak m}}		%prop.tex
\newcommand{\lz}{{{\Ell}_{\Bbb{Z}}}}	%prop.tex
\newcommand{\zq}{{\Bbb{Z}[q]}}	%prop.tex
\newcommand{\vsmt}{{v_\sigma^{**}}} 	%theorem.tex
\newcommand{\tab}{{\cal T}}		%theorem.tex
\newcommand{\us}{{\underline{s}}}	%theorem.tex
\newcommand{\dsmb}{{\frak s}}		%theorem.tex

\newcommand{\mata}{{\underline{\underline{A}}}}  %example.tex
\newcommand{\ua}{{\underline{a}}}	%example.tex
\newcommand{\uaprime}{{\underline{a'}}}	%example.tex
\newcommand{\qbinom}[2]{\left[\begin{array}{cc}{#1}\\{#2}\end{array}\right]} %example.tex
\renewcommand{\pmatrix}[1]{\left(\begin{array}{#1}}	%example.tex
\renewcommand{\endpmatrix}{\end{array}\right)}		%example.tex

\addtolength{\itemsep}{-100mm}
%\addtolength{\topsep}{-20pt}
%\input{title.tex}
%%%%%%%%%%%%%%%%%%%%%%%%%%%%%%%%%%%%%%%%%%%%%%%%%%
\title{Dual Standard Monomial Theoretic Basis and\\ Canonical Basis for Type $A$}
\author{K.N.~Raghavan %\footnotemark
  and
P.~Sankaran\footnote{Both authors were partially supported by DST under
grant \#~MS/I--73/97.  New address for P.~Sankaran:
Institute of Mathematical Sciences,  CIT Campus, Chennai 600~113, India.}
\\
Chennai Mathematical Institute \\
92 G.N.Chetty Road, T.Nagar \\
Chennai 600 017 INDIA \\
E-mail: knr@cmi.ac.in and sankaran@imsc.ernet.in
}
\date{}

%%%%%%%%%%%%%%%%%%%%%%%%%%%%%%%%%%%%%%%%%%%%%%%%%%%
\maketitle

%%%%%%%%%%%%%%%%%%%%%%%%%%%%%%%%%%%%%%%%%%%%%%%%%%

\newtheorem{theorem}{Theorem}[section]
\newcommand{\bthm}{\begin{theorem}}
\newcommand{\ethm}{\end{theorem}}

\newtheorem{proposition}[theorem]{Proposition}
\newcommand{\bpr}{\begin{proposition}}
\newcommand{\bprop}{\begin{proposition}}
\newcommand{\epr}{\end{proposition}}
\newcommand{\eprop}{\end{proposition}}

\newtheorem{definition}[theorem]{Definition}
\newcommand{\bdefn}{\begin{definition}\begin{rm}}
\newcommand{\edefn}{\end{rm}\end{definition}}

\newtheorem{example}[theorem]{Example}
\newcommand{\bexample}{\begin{example}\begin{rm}}
\newcommand{\eexample}{\end{rm}\end{example}}

\newtheorem{remark}[theorem]{Remark}
\newcommand{\bremark}{\begin{remark}\begin{rm}}
\newcommand{\eremark}{\end{rm}\end{remark}}

\newtheorem{corollary}[theorem]{Corollary}
\newcommand{\bcor}{\begin{corollary}}
\newcommand{\ecor}{\end{corollary}}

\newtheorem{exercise}[theorem]{Exercise}
\newcommand{\bex}{\begin{exercise}\begin{rm}}
\newcommand{\eex}{\end{rm}\end{exercise}}

\newtheorem{lemma}[theorem]{Lemma}
\newcommand{\blem}{\begin{lemma}}
\newcommand{\blemma}{\begin{lemma}}
\newcommand{\elemma}{\end{lemma}}
\newcommand{\elem}{\end{lemma}}

\newenvironment{proof}[1][]{{\sc Proof#1}:\ }{\ $\Box$}
\newcommand{\bpf}{\begin{proof}}
\newcommand{\bproof}{\begin{proof}}
\newcommand{\epf}{\end{proof}}
\newcommand{\eproof}{\end{proof}}

\newenvironment{prooftwo}[1][]{{\sc Proof#1}:\ }{}
\newcommand{\bpftwo}{\begin{prooftwo}}
\newcommand{\bprooftwo}{\begin{prooftwo}}
\newcommand{\epftwo}{\end{prooftwo}}
\newcommand{\eprooftwo}{\end{prooftwo}}

\newenvironment{solution}{{\sc Solution}:\ }{\ $\Box$}
\newcommand{\bsol}{\begin{solution}}
\newcommand{\esol}{\end{solution}}

\input{amssym.def}
\section{Introduction}\label{sintro}
Let $\uqminus$ be the negative part of the quantized universal
enveloping algebra constructed from a Cartan matrix associated
to a complex semisimple Lie algebra
$\frak g$.     Let $\lambda$
be a dominant integral weight and $V(\lamdba)$ the irreducible
$\uqminus$-module with highest weight $\lamdba$.   There is on the one hand
the canonical basis for $\uqminus$ \cite{kash:bc, lusz:jams, lusz:iqg}
and on the other the standard
monomial theoretic basis for the dual of $V(\lambda)$ \cite{gmodp4, gmodp5,
litt:jams}.
It is natural to ask if there is any relationship between these two bases.
To quote Littelmann \cite[page~552]{litt:jams},   ``~\ldots the properties of
the path basis suggest that the transformation  matrix should be
upper triangular \ldots''.    It is the purpose of this note to prove
that such is indeed the case when the Cartan matrix is 
of type $A$. %---see Theorem~\ref{tmain}.
 As to other types,
we have nothing to say.

Let us indicate a little more precisely what is proved here.
We show first of all that the duals of the standard monomial theoretic
bases for various $V(\lambda)$ patch together to give what can be called
a {\em dual standard monomial theoretic basis for $\uqminus$}.     
This basis lives in the crystal lattice  and 
the image modulo $q$
of a basis element is---as is perhaps to be expected---the
corresponding standard tableau thought of as a crystal.  
The main result is that the
transformation matrix between this basis and the canonical basis is
unipotent upper triangular with respect to a natural partial order on the
set of standard tableaux.      And, finally, all this holds over the 
integral forms.

These results are stated and proved in \S\ref{stheorem}.
The key to the results is the proposition proved in \S\ref{sprop}.
In \S\ref{smonom} we give a procedure to associate monomials to
tableaux on which everything else is based.   The combinatorial
properties of this procedure are stated in the lemma of \S\ref{scomb}.
These properties are crucial for the proofs.
Finally, in \S\ref{sexample}, we compute explicitly the dual standard
monomial theoretic basis for $U_q^-(sl_3)$.

We assume throughout that $\fg=sl_n(\Bbb{C})$.  We set $\ell:=n-1$
and denote by $\alpha_i$ the simple root $\epsilon_i-\epsilon_{i+1}$.
The terminology and notation of \cite{jant:lqg} are in force throughout
but for one or two minor changes in notation which should cause no confusion.

We now pass some bibliographical remarks:
\vspace{-3mm}
\begin{itemize}
\item
The most general version of standard monomial theory
is that given by Littelmann in \cite{litt:jams}.
For a readable and up to date account of standard monomial
theory and its applications, see \cite{llm:smtaa}.
Our reference for material on quantum groups and canonical basis
is \cite{jant:lqg}.
\item
As pointed out to us by Littelmann,  the procedure in \S\ref{smonom}
of associating 
a monomial to a standard tableau
is a special case of that in \cite{rs:tg}.
The crucial property of this association in the special case is
that the associated monomial is an adapted string in the 
sense of \cite{litt:ccp}.
\item
Statements (1) and (3) of Corollary~\ref{cprop} have been proved
by Littelmann \cite[Theorems~25,~17]{litt:jpaa}, at least in the
special case $q=1$.  Statement~(3)  can easily be deduced from the
results of
Berenstein and Zelevinsky~\cite{bz:string} or from those of 
Chari and Xi~\cite{cxi:preprint}.
It is also a special case of a result of Reineke \cite[see \S8]{reineke:feigin}.
Lakshmibai~\cite{lakshmi:banff} has constructed monomial bases in a very general
set up of which
statement~(3) is a special case (in this connection, see also
\cite[\S10]{litt:ccp}).  Our approach is quite different from
those of the above papers.
\item
The proof of Theorem~\ref{tsmtcan} is modelled after the
proof of Theorem~2 in \cite{cxi:preprint}.   It is noteworthy that monomial bases
play an important role in that paper as well as in Reineke's paper \cite{reineke:feigin}.
\item
It might be of interest to know how the dual standard monomial theoretic basis
relates to such other bases as the PBW basis and Reineke's dual monomial
basis \cite{reineke:feigin}.
\end{itemize}
\noindent
{\sc Acknowledgements:} \ 
We have benefitted from comments 
of Peter Littelmann on a preliminary version of this paper and
from the lectures of Jens Jantzen on quantum groups at Chennai
Mathematical Institute.   It is a 
pleasure to thank them both.   The first named author thanks
Universit\'e Louis Pasteur et Institut Universitaire de France
and the second named author thanks IHES for hospitality on visits
during which parts of this work were done.    The first named
author also thanks CNRS and NBHM for supporting his visit.
\section{Monomials and Tableaux}\label{smonom}
We will be considering monomials of a particular form in operators
indexed by the simple roots.   Let $\alpha_i$ be the simple
root $\epsilon_i-\epsilon_{i+1}$.   Consider the monomial
\[
\alpha_1^{a^1_1}\,(\alpha_2^{a^2_2}\alpha_1^{a_1^2})\,\cdots\,
	(\alpha_\ell^{a^\ell_\ell}\cdots\alpha_1^{a^\ell_1})
			\]
where $\ul{a}=(a_k^r\,|\, 1\leq r\leq \ell,\, 1\leq k\leq r)$
is any collection on non-negative integers.   Such a monomial is
{\em standard} if, for each $r$,  $1\leq r\leq \ell$,
\[
	a_r^r\geq a_{r-1}^r\geq \ldots\geq a_1^r
						\]

We now describe procedures for associating
to a tableau a standard monomial and to a standard monomial
an equivalence class of special kind of standard tableaux.
Let us first
recall the notions of tableau and standard tableau.

  Let $\shape=(m_1,\ldots, m_\ell)$
be an $\ell$-tuple of non-negative integers.   To $\lambda$ we 
associate a shape as follows.    The shape consists of boxes
$m_1+2m_2+3m_3+\cdots+\ell m_\ell$ in number,   top-justified and
right-justified,  with 1 box each in the first $m_1$ columns,  
$2$ boxes each in the next $m_2$ columns,  and so on.    For example,
if $\ell=3$,   the shape corresponding to $\shape=(3,4,2)$ is shown
in Figure~1.
\begin{center}
\begin{tabular}{|c|c|c|c|c|c|c|c|c|} 
\hline
~& ~& ~& ~& ~& ~& ~& ~& \\ \hline
\multicolumn{1}{c}{~}& \multicolumn{1}{c}{~}
& \multicolumn{1}{c|}{~}
 & ~&  ~& ~& ~& ~& \\ \cline{4-9}
\multicolumn{1}{c}{~} & \multicolumn{1}{c}{~}
& \multicolumn{1}{c}{~} & \multicolumn{1}{c}{~}
& \multicolumn{1}{c}{~} & \multicolumn{1}{c}{~}
& \multicolumn{1}{c|}{~} 
 & ~& \\ \cline{8-9}
\end{tabular}\\ [6mm]
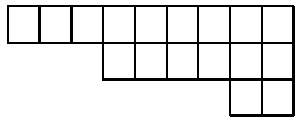
{Figure 1: the shape $(3,4,2)$.}
\end{center}

A {\em tableau} of shape $\shape$ is a filling up of  the boxes
in the shape associated to $\shape$ by integers $1,2,\ldots,
n=\ell+1$,   such that the entries in each column are strictly increasing
downwards.   A tableau is {\em standard} if the numbers in each row are
non-increasing rightwards.    For example,  if $\ell=3$ and
$\shape=(3,4,2)$,  the tableau in Figure~2 is not standard while the one
in Figure~3 is standard.
\begin{center}
\begin{tabular}{cc}
\begin{tabular}[t]{|c|c|c|c|c|c|c|c|c|} 
\hline
4& 4& 3& 3& 3& 2& 2& 2& 1\\ \hline
\multicolumn{1}{c}{~}& \multicolumn{1}{c}{~}
& \multicolumn{1}{c|}{~}
 & 4&  4& 3& 4& 3& 2\\ \cline{4-9}
\multicolumn{1}{c}{~} & \multicolumn{1}{c}{~}
& \multicolumn{1}{c}{~} & \multicolumn{1}{c}{~}
& \multicolumn{1}{c}{~} & \multicolumn{1}{c}{~}
& \multicolumn{1}{c|}{~} 
 & 4& 4 \\ \cline{8-9}
\end{tabular} \hspace{.2cm}& \hspace{.2cm} 
\begin{tabular}[t]{|c|c|c|c|c|c|c|c|c|} 
\hline
4& 4& 3& 3& 3& 2& 2& 2& 1\\ \hline
\multicolumn{1}{c}{~}& \multicolumn{1}{c}{~}
& \multicolumn{1}{c|}{~}
 & 4&  4& 4& 3& 3& 3\\ \cline{4-9}
\multicolumn{1}{c}{~} & \multicolumn{1}{c}{~}
& \multicolumn{1}{c}{~} & \multicolumn{1}{c}{~}
& \multicolumn{1}{c}{~} & \multicolumn{1}{c}{~}
& \multicolumn{1}{c|}{~} 
 & 4& 4 \\ \cline{8-9}
\end{tabular} \\ [6mm]
& \\
{Figure 2: a non-standard tableau} & {Figure 3: a standard tableau}
\end{tabular}
\end{center}

There clearly exists a smallest tableau of a given shape,
namely the one whose entries on row $r$ are all equal to $r$
for every $r$.
  The
smallest tableau of shape $\shape=(3,4,2)$ is shown in Figure~4.
\vspace{4mm}
\begin{center}
\begin{tabular}{|c|c|c|c|c|c|c|c|c|} 
\hline
1& 1& 1& 1& 1& 1& 1& 1& 1\\ \hline
\multicolumn{1}{c}{~}& \multicolumn{1}{c}{~}
& \multicolumn{1}{c|}{~}
 & 2&  2& 2& 2& 2& 2\\ \cline{4-9}
\multicolumn{1}{c}{~} & \multicolumn{1}{c}{~}
& \multicolumn{1}{c}{~} & \multicolumn{1}{c}{~}
& \multicolumn{1}{c}{~} & \multicolumn{1}{c}{~}
& \multicolumn{1}{c|}{~} 
 & 3& 3 \\ \cline{8-9}
\end{tabular} \\ [6mm]\nopagebreak
{Figure 4:  the smallest tableau of shape $(3,4,2)$}
\end{center}

Let $\sigma$ be a %standard
tableau.
For integers $r$ and $k$ such that $1\leq r\leq \ell$ and
$1\leq k\leq r$,  let $a_k^r(\sigma)$ or simply $a_k^r$ be the number
of entries
on the top
$k$ rows of $\sigma$ that are equal to $r+1$.
We define the monomial
$\monom(\sigma)$  associated to $\sigma$ to be 
\[
\monom(\sigma):= \alpha_1^{a_1^1} \left( \alpha_2^{a_2^2}\alpha_1^{a_1^2}
\right)\left( \alpha_3^{a_3^3}\alpha_2^{a_2^3}\alpha_1^{a_1^3}\right)
\cdots
\left(\alpha_\ell^{a^\ell_\ell}\cdots\alpha_1^{a^\ell_1}\right)
\]
%where $\alpha_i$ is the simple root $\epsilon_i-\epsilon_{i+1}$.
Clearly $\monom(\sigma)$ is standard.
For example, the monomial associated to the standard tableau of
Figure~3 is
\[
\alpha_1^3\,(\alpha_2^6\alpha_1^3)\,(\alpha_3^7\alpha_2^5\alpha_1^2)
	\]

We now want to associate to a standard monomial $\ul{a}=(a_k^r)$ a
standard tableau $\sigma(\ul{a})$.  Set $m_1:=a^\ell_1+\cdots+a^1_1$,
$m_j:=(a^\ell_j-a^\ell_{j-1})+\cdots+(a_j^j-a^j_{j-1})$,
$\lambda:=(m_1,\ldots,m_\ell)$,  and let $\sigma(\ul{a})$ be the standard
tableau with exactly $a^r_k-a^r_{k-1}$ entries equal to $r+1$ on row $k$.
For example,  the tableau associated to the monomial associated to the
standard tableau of Figure~3 is shown below in Figure~5.  
\begin{figure}
\begin{center}
\begin{tabular}[t]{|c|c|c|c|c|c|c|c|c|c|c|c|c|c|c|c|} 
\hline
4& 4& 3& 3& 3& 2& 2& 2& 1&1&1&1&1&1&1&1\\ 
\hline
\multicolumn{8}{c|}{~}&4&4&4&3&3&3&2&2\\
\cline{9-16}
\multicolumn{14}{c|}{~}&4&4\\
\cline{15-16}
\end{tabular} \\ [6mm]  \nopagebreak
{Figure 5:  the tableau associated to $\alpha_1^3(\alpha_2^6
\alpha_1^3)(\alpha_3^7\alpha_2^5\alpha_1^2)$}
\end{center}
\end{figure}

A tableau is {\em special} if, whenever an entry
on row $r$ is greater than $r$,  that entry is the last one in its
column.   Note that a special standard tableau  
$\sigma$ remains standard after inserting into it
the smallest tableau of shape consisting of a column of $r$ boxes
between columns $m_1+\cdots+m_r$ and $m_1+\cdots+m_r+1$.  Two special
tableaux are {\em equivalent} if one can be obtained from the other
by inserting and deleting as above smallest single column tableaux.

The {\em weight} of a tableau $\sigma$
is the element
$\sum_{r,k} a_k^r(\sigma)\alpha_k$ of the positive root lattice.
\bremark
\begin{enumerate}
\item
Standard monomials and equivalence classes of special standard tableaux
are in bijection via the maps $\sigma\mapsto\monom(\sigma)$ and
$\ul{a}\mapsto\sigma(\ul{a})$.
\item
Let $\mu$ be an element of the positive root lattice.
For every sufficiently large dominant integral weight $\lambda$,
the standard tableaux of shape $\lambda$ and weight $\mu$ are all
special.   In particular,   standard monomials of weight $\mu$
are in bijection with standard tableaux of shape $\lambda$ and
weight $\mu$ for $\lambda\gg\mu$.
\end{enumerate}
\eremark
\bproof
The first assertion is easily verified.
So is the second: if $\mu=b_1\alpha_1+\cdots+b_\ell\alpha_\ell$,
then the assertion holds for every $\lambda=(m_1,\ldots,m_\ell)$
with $m_j\geq b_j$.\nolinebreak\eproof
\section{A Combinatorial Lemma}\label{scomb}
The combinatorial properties of the association of monomials to
standard tableaux play a crucial role in the proof of the key 
Proposition~\ref{plin}.
These properties are stated in Lemma~\ref{pcomb} below,
to state and prove which is the purpose of this section.
It is convenient for this purpose
to give an alternative construction of the monomial.

The monomial associated to $\sigma$ can also be defined inductively
as follows.   Let $c$ be the least natural number such that $c+1$
occurs as an entry in some row $r$ of $\sigma$ with $r\leq c$.  In
other words,  the corresponding entry in the smallest tableau is at
most $c$.   Let us call such an entry {\em marked}.    A column
carrying a marked entry is also called {\em marked}.   By the minimality
of $c$,   the entry just above a marked entry is at most $c-1$. 
Thus, if we change all marked entries from $c+1$ to $c$,  the result
will be a tableau---let us call it $\tau$---which is also standard.
Set
\[
	\monom(\sigma):= \alpha_c^k \monom(\tau)
\]
where $k$ is the number of marked entries of $\sigma$, and $\monom(\tau)$
is defined by induction: $\tau$  is a ``smaller'' tableau---the sum of the
entries, for example,  goes down on passage from $\sigma$ to $\tau$.
The monomial associated to the smallest tableau is by definition $1$.

An entry (respectively column) of $\tau$ is called {\em marked (relative to
$\sigma$)}\  if the
corresponding
entry (respectively column) of
$\sigma$ is marked.    

Introduce a partial order on the set of all tableaux of a given
shape as follows.   For a tableau $y$, denote by $y_j$ the column $j$
 of $y$, and by $y_j(r)$ the entry on row 
$r$ of $y_j$.    If $x$ and $y$ are tableaux of shape 
consisting of a single column,   say $x\leq y$ if 
$x(r)\leq y(r)$ for every $r$.   For tableaux of general shape,
say $x\leq y$ if $x_j\lneq y_j$ for the least $j$ such that
$x_j\neq y_j$.

A tableau $y$ of shape a single column is of {\em type I} if $c$
appears as an entry in $y$ but not $c+1$,  of {\em type II} if 
either both $c$ and $c+1$ or neither appear,   and of {\em type III}
if $c+1$ appears but not $c$.   (These types correspond respectively
to the $\alpha_c$-weight being $1$, $0$, or $-1$.   Since every
fundamental weight is miniscule,  these possibilities are exhaustive.)
\blemma\label{pcomb}\label{lcomb}
Let $\sigma$ be a standard tableau.    Let $c$, $k$, and $\tau$ be as in the
inductive definition above of $\monom(\sigma)$.   Then
\begin{enumerate}
\item[(A)]
If $c+1$ occurs in $\tau$, it is only on row $c+1$. In particular,
no column of $\tau$ is of type~III.
\item[(B)]
If a column of $\tau$ containing $c$ as an entry is to the left of a marked 
column,  then that column is of type I and is itself marked. 
\item[(C)]
Let $y$ be a tableau such that $y\lneq \tau$ and $\weight(y)=
\weight(\tau)$.    Let $k$ be the number of marked entries in $\tau$.
Suppose that $x$ is a tableau obtained from $y$ as follows:
first alter $k'$  columns of $y$ of type III by replacing
$c+1$ by $c$,  where $k'$ is any non-negative integer;  then
change $k+k'$ columns 
of type I of the resulting tableau by replacing $c$ by $c+1$.
Then $x\lneq \sigma$.
\end{enumerate}
\elemma
\bproof
(A):   The first statement follows from construction---all 
$c+1$ in rows $1$ through $c$ are changed to $c$ on passage from
$\sigma$ to $\tau$.    As for the second statement, note that if
$c+1$ occurs and $c$ does not in a column,  then that $c+1$
must occur on row $i$ for $i\leq c$, which contradicts the first
statement.

(B):  
Now suppose
that $\tau_j(i)=c$ and that $\tau_j$ is not marked.   By the minimality
of $c$,  we have $i\geq c$ and so, by the tableauness of $\tau$, we get $i=c$.
By the standardness of $\sigma$,
no entry of $\sigma$ to the ``northeast'' of $\sigma_j(i)=c$  can equal
$c+1$.  This means that no column  to the right of column $j$ is marked.
If $\tau_j$ is of type II, then $\tau_j(i+1)=c+1$, and we get
$i+1\geq c+1$ just as before.    This means that in $\sigma_p(r)=r$ for $1\leq r \leq c+1$
and $p\geq j$, so that  no such column $p$ is marked.

(C):  Let $r$ be the least integer such that $y_r\neq \tau_r$.
Since $y\lneq\tau$ by hypothesis, we have $y_r\lneq\tau_r$.
Let $s$ be the least natural number such that $c$ occurs in
$\tau_s$ but $\tau_s$ is not marked.   We have two cases.

Case~1:  Assume that $r<s$.    Suppose that, for some
$j<r$,   $c$ occurs in $y_j$ and also in $x_j$.   For the
least such $j$, we clearly have $x_i=\sigma_i$ for $i<j$
and $x_j\lneq \sigma_j$,  so we are done.  We may therefore assume that,
for $j<r$,  if $c$ occurs in $y_j$,  then it changes to $c+1$
in $x_j$.   We then have $x_j=\sigma_j$ for $j<r$.   We claim
that $x_r\lneq \sigma_r$.   To prove the claim,  we may assume that
$y_r$ is of type I and that the $c$ in $y_r$ changes to $c+1$
in $x_r$, for otherwise $x_r\leq y_r\lneq\tau_r\leq \sigma_r$.
Suppose that $y_r(i)=c$.   Then clearly $i\leq c$.
It follows from (A) that either  $\tau_r(i)\geq c+2$ or $\tau_r(i)=c$.
In the former case, we clearly have $x_r(j)=y_r(j)\leq \tau_r(j)
\leq \sigma_r(j)$ for $j\neq i$, and $x_r(i)=c+1\lneq\tau_r(i)
\leq \sigma_r(I)$, so we are done.   In the latter case, we have
$x_r(j)=y_r(j)\leq \tau_r(j)=\sigma_r(j)$ for $j\neq i$, with
strict inequality holding for some $j\neq i$ (since $y_r\lneq\tau_r$
by hypothesis),  and $x_r(i)=c+1=\sigma_r(i)$,  where the last
equality holds  since $r<s$.

Case~2:   Suppose that $r\geq s$.   For $j<s$,  any $c$ occurring
in $\tau_j$ changes to $c+1$ on passage to $\sigma_j$,  so that
we have $x_j\leq \sigma_j$.   If any such $c$ does not change
on passage from $y_j$ to $x_j$,  we have $x_j\lneq\sigma_j$, and
we are done.   So we may assume that all such $c$ do change to
$c+1$ in $x$,  which means that $x_j=\sigma_j$ for $j<s$.

The $c$ that occurs in $\tau_s$ remains as such in $\sigma_s$.
By the choice of $c$,  we conclude that this $c$ occurs on
row $c$.  Thus the first $c$ rows of $\tau_j$ for $j\geq s$ are
all like those of the smallest tableau.    Since $\weight(y)=\weight(\tau)$
by hypothesis,   and $y_j=\tau_j$ for $j<s$,   it follows that the first $c$
rows of $y_j$ for $j\geq s$ are also like those of the smallest tableau.
Combining this with (A),  we conclude that $y$ does not have any
columns of type III.   So $k'=0$.   Since $y_j=\tau_j$ and $x_j=\sigma_j$
for $j<s$,  it follows that $k$ changes occur in columns 
$1$ through $s-1$ on passage from $y$ to $x$.   Thus $x_j=y_j$
for $j>s$.   In particular,  $x_j=\sigma_j$ for $j<r$ and
$x_r=y_r\lneq \tau_r=\sigma_r.$\eproof
\section{Monomial Bases}\label{sprop}
The purpose of this section is to prove Proposition~\ref{plin} below,
which provides the key to the results of \S\ref{stheorem}.
  Corollary~\ref{cprop}
provides the justification for the title of this section.

Denote by $\uqminus$ the negative part of the quantized enveloping
algebra of $\fg=sl_n(\Bbb{C})$,  by $\uzminus$ Lusztig's integral
form of $\uqminus$,   by $\linfty$ the crystal lattice of $\uqminus$,
by $A$ the local ring of fractions $f/g$ with $f$ and $g$ in the
polynomial ring $\Bbb{Q}[q]$ and $g(0)\neq 0$,  and by $\maxi$ the
maximal ideal of $A$.

Denote by $\varpi_i$ the fundamental weight $\epsilon_1+\cdots+
\epsilon_i$,  by $V(i)$ the fundamental representation associated to
$\varpi_i$,  by $v_i$ the highest weight vector $e_1\wedge\cdots
\wedge e_i$ of $V(i)$ (where $e_1,\ldots, e_n$ is the standard basis of
the standard representation $V(1)$), 
 by $\vz(i)$ the integral
form of $V(i)$ determined by $v_i$,   
 by $\Ell(i)$ the crystal lattice
of $V(i)$ determined by $v_i$,   and 
 by $\lz(i)$ the $\zq$-form $\vz(i)\cap\Ell(i)$ for $\Ell(i)$.

Let $\lambda=m_1\varpi_1+\cdots+m_\ell\varpi_\ell$ be a dominant
integral weight.  Set
\[
V:= \underbrace{V(1)\tensor\cdots\tensor V(1)}_{\mbox{$m_1$ times}}
	\,\tensor\cdots\tensor\,\underbrace{V(\ell)\tensor\cdots
		\tensor V(\ell)}_{\mbox{$m_\ell$ times}}
	\]
and define similarly $\vz$, $\Ell$, and $\lz$.  Set
\[
v_\lambda:= \underbrace{v_1\tensor\cdots\tensor v_1}_{\mbox{$m_1$ times}}
	\,\tensor\cdots\tensor\,\underbrace{v_\ell\tensor\cdots
		\tensor v_\ell}_{\mbox{$m_\ell$ times}}
	\]
Denote by $\vz(\lambda)$ the integral form of $V(\lambda)$ determined
by $v_\lambda$.   We define similarly $\Ell(\lambda)$ and $\lz(\lambda)$.
	%where $v_j$ is the highest weight vector in $\vz(j)$.

Standard tableaux of shape consisting of a single column of 
$j$ boxes index a basis for $\latz(j)$:
if $i_1<\ldots<i_j$
are the entries of a tableau $x$,  the corresponding basis element is 
$v_x := e_{i_1}\wedge\cdots\wedge e_{i_j}$.   
It follows that tableaux of shape $\lambda$ form a basis for $\lz(\lambda)$:
if $x_1,\ldots,x_m$ are the columns of a tableau $x$ of shape $\lambda$,  where
$m:=m_1+\cdots+m_\ell$,    the corresponding basis element is
$v_x:=v_{x_1}\tensor\cdots\tensor v_{x_m}$.   

Denote by $F_\alpha$ the generator of $U_q^-(\frg)$ indexed
by the simple root $\alpha$ (see \cite[4.3 and 4.4]{jant:lqg}).
The symbol $\widetilde{F}_\alpha$ will denote, depending upon
the context,   either the operator defined on a finite dimensional
$U_q^-(\frg)$-module as in \cite[9.2]{jant:lqg},  or its 
``global version'' the operator defined as in \cite[10.2]{jant:lqg}.

For a standard tableau $\sigma$,   denote by $\widetilde{F}(\sigma)$
the monomial $\monom(\sigma)$ in the operators $\widetilde{F}_\alpha$:
for instance,  if $\sigma$ is the tableau of Figure~3, we have
$\widetilde{F}(\sigma):= 
\widetilde{F}_1^3
(\widetilde{F}_2^6
\widetilde{F}_1^3)
(\widetilde{F}_3^7
\widetilde{F}_2^5
\widetilde{F}_1^2)$, where $\widetilde{F}_i$ stands for 
$\widetilde{F}_{\alpha_i}$.     Similarly,   $F(\sigma)$ denotes the
monomial $\monom(\sigma)$ in divided powers of $F_\alpha$:  for the
tableau $\sigma$ of Figure~3,  we have $F(\sigma):=
F_1^{(3)}
(F_2^{(6)}
F_1^{(3)})
(F_3^{(7)}
F_2^{(5)}
F_1^{(2)})$, where $F_i$ stands for $F_{\alpha_i}$.   
 We will now prove
that, for  a standard tableau $\sigma$ of shape $\lambda$,
the expressions for $F(\sigma)v_\lambda$ and $\widetilde{F}(\sigma)v_\lambda$
as linear combinations of the basis elements $v_x$ have a certain nice
form.

\bprop\label{plin}
For a standard tableau $\sigma$ of shape $\lambda$, we have
\begin{eqnarray}
F(\sigma)v_{\lambda}&=&\ds v_{\sigma}+\sum_{x\lneq\sigma}
                     n_{x}(\sigma)v_{x}\quad\quad\mbox{ with
$n_{x}(\sigma)\in \Bbb{N}[q,q^{-1}]$}\label{eone}\\
\widetilde F(\sigma)v_{\lambda}&=&v_{\sigma}+\sum_{x\lneq\sigma}
                       p_{x}(\sigma)v_{x}\quad\quad\mbox{ with
$p_{x}(\sigma)\in \maxi$}\label{etwo}
\end{eqnarray}
	%(Recall that $\maxi$ denotes the maximal ideal of the local ring
	%$A$ obtained by localizing $\Bbb{Q}[q]$ at the prime ideal $(q)$.)
\eprop
\bpf
If $\sigma$ is the smallest tableau,  then $F(\sigma)v_\lambda
=\widetilde{F}(\sigma)v_\lambda =v_\lambda$,  so that the statements
hold trivially.    Suppose that $\sigma$ is not the smallest tableau.
Let $c$ and $\tau$ be as in the definition of $M(\sigma)$
given in \S\ref{scomb}.    Since $\tau$ is a smaller tableau, we may
assume by induction that the statements hold for $\tau$:
\begin{eqnarray}
F(\tau)v_{\lambda}&=&\ds v_{\tau}+\sum_{y\lneq\tau}
                     n_{y}(\tau)v_{y}\quad\quad
\mbox{ with $n_{y}(\tau)\in \Bbb{N}[q,q^{-1}]$}\\
\widetilde F(\tau)v_{\lambda}&=&v_{\tau}+\sum_{y\lneq\tau}
                       p_{y}(\tau)v_{y}\quad\quad\mbox{ with
$p_{y}(\tau)\in \maxi$\label{efour}}
\end{eqnarray}
Setting $\alpha:=\alpha_{c}$ and $k:=a_{c}^{c}(\sigma)$, we have
\[
F(\sigma)v_{\lambda} = F_{\alpha}^{(k)}(F(\tau)v_{\lambda})=
                     \ds F_{\alpha}^{(k)}v_{\tau}+\sum_{y\lneq\tau}
                      n_{y}(\tau)F_{\alpha}^{(k)}v_{y}
			\]
and a similar expression for $\widetilde F(\sigma)v_{\lambda}$.

We now investigate the form of $F_{\alpha}^{(k)}v_{y}$ for  
a general tableau $y$ of shape $\lambda$. The comultiplication
$\triangle^{\prime}$ %in $U_{q}^{-}$ 
acts on $F_{\alpha}$ as follows
(see \cite[9.13~(5)]{jant:lqg}):
\[\triangle^{\prime}(F_{\alpha})=F_{\alpha}\otimes 1+K_{\alpha}
                                 \otimes F_{\alpha}\]
If $\lambda$ is a fundamental weight, then $K_{\alpha}v_{y},
F_{\alpha}v_{y}$, and $\widetilde F_{\alpha}v_{y}$ can be 
described as follows:
 \begin{eqnarray*}
F_\alpha v_y = \widetilde{F}_\alpha v_y &  = &
			\left\{\begin{array}{ll}
				0 & \mbox{if $y$ is of type II or III}\\
				v_z & \mbox{ if $y$ is of type I}
				\end{array}\right. 
								\\
	K_\alpha v_y & = & \left\{\begin{array}{ll}
				qv_y &\mbox{if $y$ is of type I}\\
				v_y  &\mbox{if $y$ is of type II}\\
				q^{-1}v_y &\mbox{if $y$ is of type III}
				\end{array}
				\right. \end{eqnarray*}
where $z$ is the tableau obtained by changing $c$ to $c+1$ in $y$
(if $y$ is of type I).

Now let $y$ be any tableau of shape $\lambda =(m_{1},\ldots,m_{\ell})$.
Set $m:=m_{1}+\ldots+m_{\ell}$.
	Define
\[
\hit(y):=\{ j\mid 1\leq j\leq m,  \mbox{column $y_j$ is of type I}\}
		\]
For a subset 
$\underline{t}=\{1\leq t_{1}<\ldots<t_{k}\leq m\}$ of cardinality 
$k$ of $\hit(y)$
define the tableau
${\underline{t}}(y)$ to be the one obtained from $y$ by replacing 
$c$ by $c+1$ in all those columns $j$ of $y$ for which $j\in\underline{t}$.
If $\hit(y)$ has exactly $k$ elements and $y$ has no columns of type III, then
\(F^{k}_{\alpha}v_{y}=[k]^!v_{{\hit(y)}(y)}  \).
In the general case,
\[F_{\alpha}^{k}v_{y}=\sum_{\underline{t}}[k]^!q^{r(\underline{t})}v_
{{\underline{t}}(y)}
\quad\mbox{ so that }\quad
F_{\alpha}^{(k)}v_{y}=\sum_{\underline{t}}q^{r(\underline{t})}v_
{{\underline{t}}(y)} \]
where $r(\underline{t}):=\sum_{i=1}^{k}(k-i+1)(\phi^{i}-\epsilon^{i})$
with $\phi^i$ and $\epsilon^{i}$ being the cardinalities respectively
of $\{j\mid t_{i-1}<j<t_{i}, \mbox{$y_{j}$ is of the type I}\}$ and
$\{j\mid t_{i-1}<j<t_{i}, \mbox{$y_{j}$ is of the type III}\}$---here
$t_0:=0$.

On subsets $\ul{t}=\{1\leq t_1<\ldots<t_k\leq m\}$ of cardinality
$k$ of $\hit(y)$,  introduce the following partial
order:    $\ul{t}\leq \ul{t'}$ if $t_j\lneq t'_j$ for the least
$j$ such that $t_j\neq t'_j$.     Let $\tnot$ be the smallest
subset of cardinality $k$ of $\hit(\tau)$.   
We claim that
\begin{enumerate}
\item[(i)]
$\sigma={\tnot}(\tau)$
\item[(ii)]
$r(\tnot)=0$
\item[(iii)]
${\ult}(y) \lneq {\ultp}(y)$ for $\ult\lneq\ultp$ in 
$\hit(y)$.
\item[(iv)]
If $y\lneq\tau$ then $\ult(y)\lneq\tnot(\tau)$
for $\ult$ in $\hit(y)$.
\end{enumerate}
It is clear that (1) of the proposition follows from the
claim.   The claim follows from Lemma~\ref{lcomb}:
(i) follows from (B) of the Lemma, (ii) from (A) and (B), and 
(iv) from (C) (with $k'=0$).   Statement (iii) is evident.

Continuing with the proof of (2),
we first observe that $E_\alpha v_\tau=0$ because of (A) of
Lemma~\ref{lcomb}, so that 
\[
\widetilde{F}_\alpha^k v_\tau =
				F_\alpha^{(k)}v_\tau =
v_\sigma + \sum_{\underline{t}\in \hit(\tau),\ \underline{t}\neq\tnot}
q^{r(\underline{t})}v_{{\underline{t}}(\tau)}
							\]
Again using Lemma~\ref{lcomb}~(A),   we see that $r(\underline{t})$
are all strictly positive.    So it remains only to worry about the
non-leading terms $p_y(\tau)v_y$ occuring in Equation (\ref{efour}).
It is enough to show that 
	\begin{equation}
\label{efive} \widetilde{F}_\alpha^k  v_y =  \sum_{x\lneq\tau} b_x(y) v_x
\quad \quad\mbox{ with $b_x(y)\in \ratf$}  
					\end{equation}
for then,  since $\widetilde{F}_\alpha^k$ preserves the lattice $\Ell$,
the $b_x(y)$ will be forced to be in $A$, and so $p_y(\tau)b_x(y)$ will
be in $\maxi$.

By definition,   $\widetilde{F}_\alpha^k v_y :=\sum_{r\geq0}
F_\alpha^{(k+r)}v_{y,r}$,  where $v_y=\sum_{r\geq 0}F_\alpha^{(r)}
v_{y,r}$ is the unique expression for $v_y$ with $v_{y,r}$ being
a vector in the highest weight space of the isotypic
$U_q(sl_2(\alpha))$-component of $V$ of highest weight $r$.
We claim that the expression for $v_{y,r}$ as a linear combination
of basis elements $v_x$ involves only such $x$ as are obtained from $y$
as follows:  first change $j$ columns of $y$ of type III by replacing in each
of them $c+1$ by $c$,  where $j\geq r$ is any integer;  then,  in the
resulting tableau,   replace $c$ by $c+1$ in any $j-r$ columns of type I.
The claim follows from the observation that the $\ratf$-span of $v_x$,
as $x$ varies over all tableaux obtained from $y$ as above for various
$r$,  is a $U_q^+(sl_2(\alpha))$-module.  
Equation~(\ref{efive}) now follows from Lemma~\ref{lcomb}~(C).\epf

\bcor\label{cprop}
\begin{enumerate}
\item
The elements $F(\sigma)v_\lambda$ as $\sigma$ runs over standard
tableaux of shape $\lambda$ form a basis for $\vz(\lambda)$.
\item
The elements $\widetilde{F}(\sigma)v_\lambda$ as $\sigma$ runs over
standard tableaux of shape $\lamdba$ form a basis for $\Ell(\lambda)$.
\item
The elements $F(\underline{a})$ as $\underline{a}$ runs over standard
monomials form a basis for $\uzminus$.
\item
The elements $\widetilde{F}(\underline{a})\cdot 1$ as $\underline{a}$
runs over standard monomials form a basis  for $\linfty$.
\end{enumerate}
\ecor
\bpf
To prove (1),   since $\vz(\lambda)$ is a free direct summand of $\vz$
of rank the number of standard tableaux of shape $\lambda$,  it is enough
to show that $F(\sigma)v_\lambda$ form part of a basis of $\vz$.  But this
is immediate from (1) of the proposition.  The proof of (2) is similar.

To prove (3),   let $\mu$ be an element of the positive root lattice
(in other words,  $\mu$ is a non-negative linear combination of the
simple roots).   Choose a dominant integral weight $\lambda$ so large
that the $\uzminus$-module map $\uzminus\rightarrow \vz(\lambda)$ given
by $1\mapsto v_\lambda$ restricts to an isomorphism of the weight space
$(\uzminus)_{-\mu}$ onto $(\vz(\lambda))_{\lambda-\mu}$,  and the
standard monomials of weight $\mu$ are in bijection with standard tableaux
of shape $\lambda$ and weight $\mu$.    By (1),   the elements $F(\sigma)
v_\lambda$ as $\sigma$ varies over standard tableaux of shape $\lambda$
and weight $\mu$ form a basis for $(\vz(\lambda))_{\lambda-\mu}$, so we
are done.

To prove (4),   we reduce as in the proof of (3) to showing that 
$(\widetilde{F}(\sigma)\cdot 1)v_\lambda$ as $\sigma$ varies over
standard tableaux of shape $\lambda$ and weight $\mu$ form a basis for
$(\Ell(\lambda))_{\lamdba-\mu}$.   By (2),  we know that 
$\{\widetilde{F}(\sigma)v_\lambda\}$ form a basis for $(\Ell(\lambda))_{
\lambda-\mu}$.  On the other hand,   by \cite[Proposition~10.9]{jant:lqg},
$(\widetilde{F}(\sigma)\cdot1)v_\lambda=\widetilde{F}(\sigma)v_\lambda
\mbox{ mod } q\Ell(\lambda)$.   We are therefore done by applying Nakayama.\epf

\section{The Theorem}\label{stheorem}
We keep the notations of \S\ref{sprop}.
By taking the transpose of the embedding $\vz(\lambda)\hookrightarrow
\vz$ (respectively $\Ell(\lambda)\hookrightarrow\Ell$,
respectively $\latz(\lambda)\hookrightarrow\latz$),  we get a
surjective mapping $\vz^*\surjection\vz(\lambda)^*$
(respectively
$\Ell^*\surjection\Ell(\lambda)^*$,
respectively $\latz^*\surjection\latz(\lambda)^*$).
Let $\{v_x^*\}$ be the dual basis in $\latz^*$ of the basis
$\{v_x\}$ of $\latz$.     It follows from Proposition~\ref{plin}~(1)
and Corollary~\ref{cprop}~(1) that the images of $v_\sigma^*$ in
$\vz(\lambda)^*$ (by abuse of notation also denoted $v_\sigma^*
$),  as $\sigma$ varies over standard tableaux of
shape $\lambda$, form a basis for $\vz(\lambda)^*$.
Similarly it
follows from Proposition~\ref{plin}~(2)
and Corollary~\ref{cprop}~(2) that the $v_\sigma^*$ form a basis
for $\Ell(\lambda)^*$.      Thus the $v_\sigma^*$ form a basis for
$\latz(\lambda)^*$.   Consider the dual basis $v_\sigma^{**}$ in
$\latz(\lambda)$.   If the weight $\mu$ of $\sigma$ is small compared
to $\lambda$,  we can think of $v_\sigma^{**}$ as an element of
$(\latz(\infty))_{-\mu}$.   We claim that this is independent of the
choice of $\lambda$:
\bprop\label{psmt}
Let $\mu$ be an element of the positive root lattice.  Let $\lambda$
be a dominant integral weight so large that $(\latz(\infty))_{-\mu}$
can be identified with $(\latz(\lambda))_{\lambda-\mu}$,  and the
standard tableaux of shape $\lambda$ and weight $\mu$ are all special.
Let $\lambda'$ be another such weight,   and let $\sigma\leftrightarrow
\sigma'$ denote the bijective correspondence between standard tableaux
of weight $\mu$ of shape $\lambda$ on the one hand and of shape
$\lambda'$ on the other.    Then $v_\sigma^{**}=v_{\sigma'}^{**}$
as elements of $(\latz(\infty))_{-\mu}$. 
\eprop
\bpf
Evaluating both sides of Proposition~\ref{plin}~(1) on $v_\nu^*$,
as $\nu$   varies over standard tableaux of shape $\lambda$,
we find that 
\begin{eqnarray}\label{esix}
F(\sigma)=v_\sigma^{**}+\sum_{\theta\lneq\sigma}n_\theta(\sigma)v_\theta^{**}
\end{eqnarray}
in $(\uzminus)_{-\mu}$, where the sum is taken only over standard tableaux
$\theta\lneq\sigma$, and similarly
\[
F(\sigma')=v_{\sigma'}^{**}+\sum_{\theta'\lneq\sigma'}n_{\theta'}(\sigma')v_{\theta'}^{**}
\]
We have $F(\sigma)=F(\sigma')$ by hypothesis.   We will presently show that
$n_\theta(\sigma)=n_{\theta'}(\sigma')$.    It will then follow that
$\{v_\sigma^{**}\}$
and $\{v_{\sigma'}^{**}\}$ are related to the basis $\{F(\sigma)\}$ by
the same transformation matrix,   which means $v_\sigma^{**}
=v_{\sigma'}^{**}$.

The following proof that $n_\theta(\sigma)=n_{\theta'}(\sigma')$ looks more
difficult than it really is.    It is easy if one thinks in terms of
pictures,  but to express it in words requires cumbersome notation.
We may assume, without loss of generality,  that $\lambda'\geq\lambda$,
that is,  $\lambda=m_1\varpi_1+\cdots+m_\ell\varpi_\ell$ and
$\lambda'=m_1'\varpi_1+\cdots+m_\ell'\varpi_\ell$, with $m_1'\geq m_1,
\ldots, m_\ell'\geq m_\ell$.   Given a tableau $y$ of shape $\lambda$,
we associate to it a tableau $y'$ of shape $\lambda'$ as follows:
for each $r$, $1\leq r\leq\ell$, insert into $y$,  between columns
$m_1+\cdots+m_r$ and $m_1+\cdots+m_r+1$,  $m_r'-m_r$ columns each equal
to the smallest tableau of shape consisting of a column of $r$ boxes.
The association $y\rightarrow y'$ is injective,  it generalizes the
association $\sigma\leftrightarrow\sigma'$,  and it preserves the 
property of being special.  Denote by $\tab$ the set of tableaux of
shape $\lambda'$ that are obtained as $y'$ from special $y$ of shape
$\lambda$.    Set $n_\sigma(\sigma):= 1 $ and $n_y(\sigma):=0$ for
$y\nleq\sigma$.

We will prove the following slightly stronger statement:
\[
	F(\sigma')v_{\lambda'} = \sum\limits_{y'\in\tab}
			n_y(\sigma)v_{y'} + \sum\limits_{w\not\in\tab}
			n_w(\sigma')v_w
					\]
Let $c$, $k$, $\alpha=\alpha_c$, and $\tau$ be as in the inductive 
definition of $\monom(\sigma)$ in \S\ref{scomb}.   We may assume,
by way of induction, that the statement holds for $\tau$:
\[
	F(\tau')v_{\lambda'} = \sum\limits_{z'\in\tab}
			n_z(\tau)v_{z'} + \sum\limits_{x\not\in\tab}
			n_x(\tau')v_x
					\]
We have, by definition,
\[
	F(\sigma')v_{\lambda'} : = F_\alpha^{(k)}(F(\tau')v_{\lambda'})
		= \sum n_z(\tau) F_\alpha^{(k)} v_{z'} + \sum
				n_x(\tau') F_\alpha^{(k)}v_x  \]
It is convenient to use again the notation introduced in the proof
of Proposition~\ref{plin}.     The following statements are evident:
\begin{itemize}
\item
for $x\not\in\tab$ and $\us\in\hit(x)$, we have $\us(x)\not\in\tab$.
\item
for $z'\in\tab$ and $\us\in\hit(z')$,  the tableau $\us(z')$ belongs
to $\tab$ if and only if $\us$ is of the form $\underline{t}'$ for
some $\ult\in\hit(z)$ such that $\ult(z)$ is special:  for $\ult=
\{1\leq t_1<\ldots<t_k\leq m\}$, we define $\ultp :=
\{1\leq t_1'<\ldots<t_k'\leq m'\}$ (where $m:= m_1+\cdots+m_\ell$
and $m':= m_1'+\cdots+m_\ell'$) by $t_j'=t_j+(m_1'-m_1)+\cdots+
(m_p'-m_p)$,  where $p$, $0\leq p\leq \ell-1$, is such that
$m_1+\cdots+m_p < t_j \leq m_1+\cdots+m_{p+1}$.
\end{itemize}

It therefore remains only to show that $r(\ult)=r(\ultp)$ for
$z'\in\tab$ and $\ult\in\hit(z)$ such that $\ult(z)$ is special.
Now, since $\ult(z)$ is special,  we have $t_k\leq m_1+\cdots + m_c$.
And, since $z'$ is in the image of the association $y\rightarrow y'$,
for any $p$, $1\leq p \leq c-1$, and any $j$ such that
$m_1'+\cdots +m'_{p-1}+m_p < j \leq m_1'+\cdots+m_p'$, the column
$z_j'$ is of type II.   It should now be clear that $r(\ult)=r(\ultp)$.\epf

\noindent
It follows from the proposition above that to each standard monomial
$\underline{a}$ we can associate an element $\dsmb(\underline{a})$
of $(\latz(\infty))_{-\mu}$: set $\dsmb(\underline{a}):=v_\sigma^{**}$,
where $\sigma$ is the standard tableau of shape $\lambda$ with associated
monomial $\underline{a}$, and $\lambda\gg\mu$.   
The elements $\dsmb(\underline{a})$,  as $\underline{a}$ varies over standard
monomials,  form a basis for $\latz(\infty)$.   We call this the
{\em dual standard monomial theoretic basis}.

We claim that the element $\vsmt$ of $\linfty$ maps to $\sigma$ modulo $q$.
To prove this,  evaluate both sides of Proposition~\ref{plin}~(2) on
$v_\nu^{*}$ as $\nu$ varies over standard tableaux of shape $\lambda$ to
get 
\[
	\widetilde{F}(\sigma)v_\lambda=\vsmt+\sum_{\tau\lneq\sigma}
				p_\tau(\sigma)v_\tau^{**}
								\]
where the sum is taken only over standard tableaux $\tau\lneq\sigma$.
Choosing $\lambda$ large compared to the weight of $\sigma$,  we may assume that
$\vsmt$ and $v_\tau^{**}$ in the last equation are the images of the corresponding elements
$\vsmt$ and $v_\tau^{**}$ in the algebra under the map $1\mapsto v_\lambda$.
Since $\widetilde{F}(\sigma)v_\lambda$ maps to $\sigma$ modulo $q$ and
$p_\tau(\sigma)$ vanish modulo $q$,  the claim follows.

It is natural to ask for the relation between the image of the algebra
element $\vsmt$ under the map $1\mapsto v_\lambda$ on the one hand and
the element $\vsmt$ of $V(\lambda)$ on the other when $\lambda$ is not
necessarily large compared to the weight of $\sigma$.   Since the algebra basis
$\{\vsmt\}$ is unipotent upper triangular related to $\{F(\sigma)\}$ and
$\{F(\sigma)v_\lambda\}$ is unipotent upper triangular related to the module basis
$\{\vsmt\}$,  it follows that the matrix relating the bases is 
unipotent upper triangular.   Furthermore, since both live in $\lz(\lambda)$,
the coefficients of the matrix are in $\zq$.    And, since both $\vsmt$ map to 
$\sigma$ modulo $q$,  the entries of this matrix strictly above the diagonal
are all divisible by $q$.
\bthm\label{tsmtcan}
The transformation matrix between the dual standard monomial theoretic
basis $\{v_\sigma^{**}\}$  and the canonical basis $\{G(\sigma)\}$
is unipotent upper triangular with respect to the partial order on the
standard tableaux defined in \S\ref{scomb}.    The entries of this
matrix are in $\zq$ and those strictly above the diagonal
are all divisible by $q$.
\ethm
\bpf
We first note that
the second assertion follows easily from the first.
Suppose that the matrix relating the two bases is unipotent upper
triangular with entries in $\ratf$.   Then, 
since both bases live in $\lz(\infty)$,   it follows that the entries of
the transformation matrix belong
to $\Bbb{Z}[q,q^{-1}]\cap A=\zq$.   Further, since both $v_\sigma^{**}$ and
$G(\sigma)$ map to $\sigma$ modulo $q$,   it follows  that the 
entries strictly above the diagonal are all divisible by $q$.

To prove the first assertion,
we concentrate on a single weight space $(\uzminus)_{-\mu}$. 
Choosing $\lambda\gg\mu$,   we may pass to $(\vz(\lambda))_{\lambda-\mu}$.
Since $\{v_\sigma^{**}\}$ is unipotent triangular related to 
$\{F(\sigma)v_\lambda\}$,  it is enough to show that $\{F(\sigma)v_\lambda\}$
and $\{G(\sigma)v_\lambda\}$ are unipotent upper triangular related.

From \cite[Proposition~10.9]{jant:lqg}, we have
\[
	\{(\widetilde{F}(\sigma)\cdot 1 )v_\lambda\} = \matp \{
		\widetilde{F}(\sigma)v_\lambda\}
	\]
where $\matp$ is a matrix with entries in $A$ and equals the identity matrix
modulo $q$.     It follows from Proposition~\ref{plin} that
\[
	\{\widetilde{F}(\sigma)v_\lambda\} = \matb \{F(\sigma)v_\lambda\}
	\]
where $\matb$ is a unipotent upper triangular matrix with entries in 
$\ratf$.     It can be proved by elementary means %(see Lemma~\ref{lelem} below)
 that such a matrix $\matb$ factorises as follows:
\[
	\matb = \matc\,\matd
				\]
where $\matc$ and 
$\matd$ are both unipotent upper triangular,   $\matc$ has entries
in $A$ and is identity modulo $q$,  and $\matd$ is {\em bar-invariant},
that is,   it does not change under the $\Bbb{Q}$-automorphism of $\ratf$
that interchanges $q$ and $q^{-1}$.

Noting that $\matp\,\matc$ is invertible (since it is so modulo $q$), we
get from these three equations
\[
	(\matp\,\matc)^{-1}\{\widetilde{F}(\sigma)\cdot1\} = \matd\{F(\sigma)\}
	\]
The left side maps to the crystal basis modulo $q$,  while the right side
is bar-invariant.   A characterisation of
canonical basis \cite[Theorem~11.10~(a)]{jant:lqg}
says that either side of the last equation is the canonical basis, so
\[
\{G(\sigma)\}=\matd\{F(\sigma)\}
				\]
Since $\matd$ is upper triangular, we are done.\epf

\section{An Example}\label{sexample}
 The purpose of this section is to calculate the dual standard monomial
theoretic basis for the Cartan matrix $A_2$,  in other words, for
$U^-_q(sl_3)$.   In this case,  the standard monomials are
\[
	\ua =\left\{(a,b,c)\st a\geq0, b\geq c\geq0\right\}
							 	\]
It follows, from Corollary~\ref{cprop}~(4) for example,  that 
\[
	\left\{F_1^{(a)} F_2^{(b)} F_1^{(c)} \st
		 a\geq0, b\geq c\geq0\right\}
							 	\]
form a basis for $U_{\Bbb{Z}}^-(sl_3)$.   
We will express the dual standard monomial theoretic basis 
\[
	\left\{\dsmb(\ua) \st \ua=(a,b,c),
		 a\geq0, b\geq c\geq0\right\}
							 	\]
in terms of this monomial basis.

Identifying the standard monomials with equivalence classes of special
standard tableaux,  and transferring to standard monomials the partial
order on tableaux defined in \S\ref{scomb},  we say
$\ua:=(a,b,c) \leq \uaprime:= (a',b',c')$ if either $c<c'$ or
($c=c'$ and $a<a'$) or ($c=c'$, $a=a'$,  and $b\leq b'$).   We need to
 make comparisons
between standard monomials $\ua$ and $\uaprime$ using this partial order only in the case
when their weights are the same,  that is,  when $b=b'$ and $a+c=a'+c'$.

Equation~(\ref{esix}) gives us 
\[
F(\ua) =\dsmb(\ua) +\sum\limits_{\uaprime<\ua} n_{\uaprime}(\ua) \dsmb(\uaprime)
									\]
where $\ua$ is a fixed standard monomial,  the sum is over all standard
monomials $\uaprime$ that have the same weight as $\ua$ and satisfy
$\uaprime<\ua$ in the above partial order,  and by $n_{\uaprime}(\ua)$
we mean $n_\theta(\sigma)$ (see Equation~(\ref{eone}))
where $\sigma$ and $\theta$ are standard
tableaux (of shape $\lambda$ large relative to the weight of $\ua$, that 
is, $\lambda=(m_1, m_2)$ with $m_1\geq a+c$ and $m_2\geq b-c$)
corresponding respectively to $\ua$ and $\uaprime$.  Recall that the
point of Proposition~\ref{psmt} is that $n_{\uaprime}(\ua)$ is 
independent of the choice of $\lambda$.      Our task then is to 
compute $n_{\uaprime}(\ua)$.   
\bprop\label{pexample}
Let $b$ and $k$ be fixed non-negative integers.    For integers
$s$ and $t$ such that $0\leq t \leq s \leq \min\{b,k\}$,  setting
$\ua=(k-s,b,s)$ and $\uaprime=(k-t,b,t)$,  we have 
\[
	n_{\uaprime}(\ua) = q^{(s-t)(b-t)}
	\qbinom{k-t}{s-t}
%\left[{k-t}\genfrac{s-t}\right]
								\]
\eprop
\bpftwo
We give a sketch of the proof---in fact, we sketch two proofs.
We have 
\[
	F(\sigma)=F_1^{(k-s)} F_2^{(b)} F_1^{(s)}
							\]
Keeping track of terms in the expansion of $F(\sigma)v_\lambda$ that can
give rise to $v_\theta$,  we get
\begin{eqnarray*}
	n_{\uaprime}(\ua)%& = & q^{(s-t)(b-2t-k-1)}
		& = & q^{(s-t)((b-t)-(k-t)-1)}
		\sum\limits_{1\leq i_1<\ldots<i_{s-t}\leq k-t}
		q^{2(i_1+\cdots+i_{s-t})}\\
\end{eqnarray*}
The proposition now follows from the following identity 
(Equation~1.3.1(c) of \cite{lusz:iqg}) which is proved
easily by induction:
\[ 		q^{s(b-k-1)}
		\sum\limits_{1\leq i_1<\ldots<i_{s}\leq k} 
		q^{2(i_1+\cdots+i_{s})} = q^{sb}\qbinom{k}{s}. \]

Actually,  there is no need to keep such careful track of
the coefficients of various relevant terms after application
of $F_1^{(s)}$.   An easier proof is obtained as follows: 
after application of $F_2^{(b)}$, each relevant term picks
up a factor $q^{(s-t)(b-t)}$,  and the quantum binomial factor
on the right side of the equation of the proposition
is accounted for by the obvious identity 
\[
F_1^{(k-s)}F_1^{(s-t)} = \qbinom{k-t}{s-t}   F_1^{(k-t)}.\quad\Box\]
\epftwo

The expression in matrix form of the proposition %(assuming $b\geq k$)
 is
\[
\mata \pmatrix{c}
	\dsmb(k,b,0)\\
	\dsmb(k-1,b,1)\\
	\vdots\\
	\dsmb(k-\min\{b,k\},b,\min\{b,k\})	\endpmatrix
	=
	\pmatrix{c}
		Y_1^{(k)}Y_2^{(b)}Y_1^{(0)}\\
		Y_1^{(k-1)}Y_2^{(b)}Y_1^{(1)}\\
		\vdots\\
		Y_1^{(k-\min\{b,k\})}Y_2^{(b)}Y_1^{(\min\{b,k\})}\endpmatrix
							\]
where 
\[
	{\mata}_{s+1,t+1}
			=\left\{ \begin{array}{ll}
			q^{(s-t)(b-t)}\qbinom{k-t}{s-t} & \mbox{if $s\geq t$}\\
			0				& \mbox{otherwise}
					\end{array}
										\right.
						\]
To obtain an expression for the dual standard monomial theoretic
basis in terms of the monomial basis,  we need to compute the inverse
of the matrix $\mata$.    We claim that the inverse $\mata^{-1}$ is given by
\[
	\mata^{-1}_{s+1,t+1} 
			=\left\{ \begin{array}{ll}
			(-1)^{s-t}q^{(s-t)(b-s+1)}\qbinom{k-t}{s-t} & \mbox{if $s\geq t$}\\
			0				& \mbox{otherwise}
							\end{array}
										\right.
						\]
To see this,  we need only verify that, for $s\geq t$,
\begin{eqnarray*}
&&	\sum_{t\leq j\leq s} \mata^{-1}_{s+1,j+1}\cdot\mata_{j+1,t+1}\\
&&	 = 	\sum_{t\leq j\leq s} (-1)^{s-j} q^{(s-j)(b-s+1)+(j-t)(b-t)}
		\qbinom{k-j}{s-j}\qbinom{k-t}{j-t}\\
&&	 = 	\sum_{0\leq j\leq s-t} (-1)^{s-j-t} q^{(s-j-t)(b-s+1)+j(b-t)}
		\qbinom{k-j-t}{s-j-t}\qbinom{k-t}{j}\\
&&	 = 	(-1)^{s-t}\qbinom{k-t}{s-t}q^{(s-t)(b-s+1)}
		\sum_{0\leq j\leq s-t} (-1)^j q^{j(s-t-1)}
		\qbinom{s-t}{j}	
		\end{eqnarray*}		
and, as can be seen by a routine induction, the sum in the last line above
is $0$ for $s>t$ and $1$ for $s=t$
\cite[Equation~1.3.4(a)]{lusz:iqg}.

\end{document}